\documentclass[11pt]{article}
\usepackage{amssymb}
\usepackage{color}
\usepackage[dvips]{graphicx}
\usepackage[dvips]{epsfig}
\usepackage{amsmath}

\newtheorem{thm}{Theorem}[section]
\newtheorem{defn}{Definition}[section]
\newtheorem{exa}{Example}

\newtheorem{prop}{Proposition}[section]
\newtheorem{cor}[thm]{Corollary}

\newtheorem{rem}{Remark}[section]
\newtheorem{lem}{Lemma}[section]

\newtheorem{crit}{Criterion}

\def\sep{\setlength{\arraycolsep}}

\def\ben{\begin{enumerate}}
\def\een{\end{enumerate}}
\def\bexa{\begin{exa}}
\def\eexa{\end{exa}}
\def\eproof{\hfill$\rule{2mm}{2mm}$\par}

\def\bqu{\begin{quote}}
\def\equ{\end{quote}}
\def\p{\partial}

\def\¯{\hbox{ƒ\it y}}
\def\{\hbox{ƒ\it z}}

\def\ds{\displaystyle}

\def\sep{\setlength{\arraycolsep}}

\def\proof{\noindent{\bf Proof: }}
\def\eproof{\hfill$\rule{2mm}{2mm}$\par}

\def\g{\hbox{$\gamma$}}
\def\d{\hbox{$\delta$}}

\def\i{\hbox{$\mathbf{i}$}}
\def\j{\hbox{$\mathbf{j}$}}
\def\k{\hbox{$\mathbf{k}$}}

\def\b{{\bf b}}

\def\be{\begin{equation}}
\def\ee{\end{equation}}
\newcommand{\ba}{\begin{eqnarray}}
\newcommand{\ea}{\end{eqnarray}}
\def\ds{\displaystyle}

\def\R{\hbox{$\mathbb{R}$}}

\def\N{\hbox{$\mathbb{N}$}}

\def\C{\hbox{$\mathbb{C}$}}
\def\H{\hbox{$\mathbb{H}$}}

\def\gcd{\hbox{\rm{gcd}}}

\def\p{{\bf p}}

\def\c{{\gamma}}
\def\d{{\delta}}

\def\b{\hbox{$\beta$}}
\def\a{\hbox{$\alpha$}}

\def\p{\partial}

\begin{document}
\baselineskip=16.5pt

\title{\Huge{\sf Cauchy-Riemann equations and Jacobians of quaternion polynomials}}

\author{\Large{Takis~Sakkalis and Sofia~Douka}
}

\date{September 17, 2016}

\maketitle

\begin{abstract}
\noindent
A map $f$ from the quaternion skew field $\H$ to itself, can also be thought as a transformation $f:\R^4 \to \R^4$. In this manuscript, the Jacobian $J(f)$  of $f$ is computed, in the case where $f$ is a quaternion polynomial. As a consequence, the Cauchy-Riemman equations for $f$ are derived. {\sf It} is also shown that the Jacobian determinant of $f$ {\sf is} non negative over $\H$. The above commensurates   well with the theory of analytic functions of one complex variable. 
{{\em MCS 2010:} 26B10, 12E15, 11R52}
\end{abstract}


\vskip 40mm
\footnoterule{\footnotesize{Takis Sakkalis, 
Mathematics Laboratory, Agricultural University of Athens, 
75 Iera Odos, Athens 11855, GREECE, corresponding author: e-mail stp@aua.gr \newline 
{\sf S}ofia Douka, 37 Strymonos St., Athens 11855, GREECE, email sofia.nefeli@gmail.com}}

\newpage

\setcounter{page}{2}

\thispagestyle{plain}

\section{Introduction} 
\label{sec:intro}

Let $\H$ denote the skew field of quaternions. Its elements are of the form $c=c_0 + \i c_1 + \j c_2 + \k c_3$, where $c_m \in \R$ and $\i, \j, \k$ are such that $\i^2=\j^2=\k^2=-1$ and $\i \j =-\j \i=\k, \j \k =-\k \j =\i, \k\i =-\i \k =\j$. 
The real part of $c$ is $Re(c)=c_0$ while the imaginary part $Im(c)=\i c_1 + \j c_2 + \k c_3$. The norm of $c$, $|c|=\sqrt{c_0^2+c_1^2+c_2^2+c_3^2}$, its conjugate $c^*=c_0 - \i c_1 - \j c_2 - \k c_3$ while its  inverse is $c^{-1}=c^* \cdot |c|^{-2}$, provided that $|c| \ne 0$. $c$ is called an {\em imaginary} unit if $Re(c)=0$ and $|c|=1$, and it has the property $c^2=-1$. In that regard, $\H$ is a real normed division (non commutative) algebra. An element $c=c_0 + \i c_1 + \j c_2 + \k c_3$ of $\H$ can also 
be represented via a real $4 \times 4$ matrix ${\cal C}$ as follows: 

$${\cal C}=\left[ \begin{array}{rrrrrr} c_0&-c_1&-c_2&-c_3\\c_1&c_0& -c_3&c_2\\c_2&c_3&c_0&-c_1\\c_3&-c_2&c_1&c_0 \end{array} \right] $$
Notice that $|{\cal C}|=|c|^4$. The following notation will be frequently used in the sequel: 

\begin{defn} For any $4 \times 4$ real matrix $B$ and a quaternion $c$, we define: $cB\equiv{\cal C} B$ and $Bc \equiv B {\cal C}$. 
\end{defn}

We may identify $\H$ with $\R^4$ via the map $(x +\i y + \j z +\k w) \to (x,y,z,w)$. Let $f:\H \to \H$. In view of this identification, we can also think of $f$ as a map from  $\R^4 \to \R^4$. Indeed, if $f(x +\i y + \j z +\k w)=f_1(x,y,z,w) + \i f_2(x,y,z,w) +\j f_3(x,y,z,w) +\k f_4(x,y,z,w)$ we define $f:\R^4 \to \R^4$ by  $f(x,y,z,w)=(f_1, f_2, f_3, f_4)$.  In that case, we have

\begin{defn} \label{jac} The Jacobian $J(f)(c),\; c \in \H$ is the matrix $\left [ \frac{\partial f_i}{\partial x_j} \right], i=1,2,3,4,$  $x_1=x, x_2=y, x_3=z, x_4=w$ evaluated at $c$. The determinant of $J(f)$ will be denoted by $|J(f)|$. 
\end{defn}

In this manuscript, we are concerned with the computation of $J(f)$ and its determinant at a point $c \in \H$, in the case where $f$ is a polynomial with coefficients in $\H$. More specifically, we will consider polynomials of the form 

 \be \label{left} 
f(t)=a_nt^n + a_{n-1}t^{n-1} + \cdots + a_0, \hbox{ where $a_k \in \H$} \ee

We will first show that $|J(t^n)| \geq 0$ at all points $c \in \H$. By extending the above, we will find the form of $J(f)$ at a complex number $r +\i s$ (the complex number $i$ is to be identified with the quaternion $\i$ throughout this note). Using this we deduce the form of $J(f)(c), \; c \in \H$. From the above, the Cauchy-Riemman equations for $f$ are derived. Finally, we prove that $|J(f)| \geq 0$ over $\H$,  a fact similar to the one in the theory of analytic functions of a complex variable. The latter can be used to compute the zeros $\zeta_k$ of $f$ as well as its local degree at $\zeta_k$.

\section{A brief overview of quaternion polynomials}

In this paragraph we will recall some facts, needed for the rest of the paper, concerning quaternion polynomials. For more details, the reader is referred to \cite{GS, gordon64, topu09}. 

Due to the non commutative nature of $\H$, polynomials over $\H$ are usually distinguished into the following three types: {\em left}, {\em right} and {\em general}, \cite{gordon64}. A left polynomial is an expression of the form (\ref{left}). If $a_n \ne 0$, $n$ is called the degree of $f$. Here we shall consider, unless otherwise stated, left polynomials only and called them simply polynomials. If $g(t)= b_mt^m + b_{m-1}t^{m-1} + \cdots + b_0$  is another polynomial, their product $fg(t)$ is defined in the usual way: 

$$ \label{prod} fg(t)=\sum_{k=0}^{m+n}c_k t^k, \quad \hbox{where $c_k=\ds\sum_{i=0}^k a_i b_{k-i}$} $$
Note that in the above setting the multiplication is performed as if the coefficients  were chosen in a commutative field.  

An equivalent representation {\sf of} $f$ is $f(t)=a(t) + \i b(t) +\j c(t) + \k d(t)=a +\i b +\j (c -\i d)$, where $a,b,c,d \in \R[t]$. If $c \in \H$, we define $f(c)=a_nc^n + a_{n-1}c^{n-1} + \cdots + a_0$; if $f(c)=0$, $c$ is called a {\em zero} or a {\em root} of $f$. According to  Theorem 1 of \cite{gordon64} an element $c\in \H$ is a zero of $f$ if and only if there exists a (left) polynomial $g(t)$ such that $f(t)=g(t)(t-c)$. 
In their nominal paper \cite{EN}, Eilenberg and Niven, using a degree argument, proved the fundamental theorem of algebra for quaternions, namely that any quaternion polynomial of positive degree $n$ of the form $f(t)=a_0ta_1t \cdots ta_n + \phi(t), a_i \in \H, a_i \ne 0$ and $\phi(t)$ is a sum of finite number of similar monomials $b_0tb_1t \cdots tb_k, k <n$, has a root in $\H$. 

Roots of $f$ are distinguished into two types: (i) {\em isolated} and (ii) {\em spherical}. A root $c$ of $f$ is called spherical if and only if its characteristic polynomial $q_c(t)=t^2-2t\, Re(c)+ |c|^2$ divides $f$; for any such polynomial, call $\a_c \pm \i \b_c$ its complex roots. In that case any quaternion $\gamma$ similar to $c$, is also a root of $f$; ($c_1, c_2 \in \H$ are called {\em similar}, and denoted by $c_1 \sim c_2$,  if $c_1 \,\eta =\eta \, c_2$ for  a non zero $\eta \in \H$). For example, if $f(t)=t^2 +1$, any imaginary unit quaternion $c=\a_1\, {\bf i} + \a_2\, {\bf j} + \a_3\, {\bf k}$  is a root of $f$. 

\begin{rem} \label{rem1} The polynomial $f(t)$ has a spherical root if and only if it has roots $\a + \i \b, \a - \i \b, \; \a, \b \in \R, \b \ne 0$. 
\end{rem}

The above facts, as well as the representation of $f$ as in (\ref{left}), allows us to factor $f$ into a product of linear factors $(t-c_i), \; c_i \in \H$. Indeed, since  $f(t)=g(t)(t-c)$ and $g(t)$ has a root, simple induction shows that 
\be \label{factor} f(t)=a_n(t-c_n)(t-c_{n-1}) \cdots (t-c_1), \; c_j \in \H \ee

A word of caution: In the above factorization, while $c_1$ is necessarily a root of $f$, $c_j, j=2, \cdots n,$ might not be roots of $f$. For example, the polynomial $f(t)=(t + \k)(t + \j)(t + \i)=t^3+ (\i +\j + \k)t^2+(-\i +\j  - \k)t +1$ has only one root, namely $t=-\i$. Theorem 2.1 of \cite{GS} provides a more detailed version of the above factorization.

If we write $f$ in the form $f(t)=a(t) + \i b(t) +\j c(t) + \k d(t)$ we see that $f$ has no spherical roots if and only if $\gcd(a,b,c,d)=1$. Such an $f$ will be called {\em primitive}. Then, it is known (Corollary 3.3, \cite{topu09}) that a primitive $f(t)$ of degree $n$, has at most $n$ distinct roots in $\H$. 

The conjugate $f^*$ of $f$ is defined as $f^*=a(t) - \i b(t) -\j c(t) -\k d(t)$. Note that $f^*  f= a^2+b^2+c^2+d^2$, which is a real positive polynomial. Observe that if $\a +\i \b$ is a root of $f^*f$, then there exists a $c \in H$, similar to $\a +\i \b$ so that $f(c)=0$.  

\begin{defn} \label{defn1} Let $\phi(t) \in \C[t]$ and $\zeta \in \C$ be a root of $\phi$. We denote by $\mu(\phi)(\zeta)$ the multiplicity of $\zeta$. Now let $c \in \H$ be a root of $f$ and let $m=\mu(f^*f)(\a_c + \i \b_c)$. Then, if (1) $c$ is isolated, we define its multiplicity $\mu(f)(c)$, as a root of $f$, to be $m$; (2) if $c$ is spherical, its multiplicity is set to be $2m$. 
\end{defn} 
Note that the above notion of multiplicity agrees with the one given in Definition 2.6 of \cite{GS}, page 23. From the above we have: 
\begin{crit} \label{crit1} Let $c$ be a root of the {\em primitive} polynomial $f$. Then $c$ has multplicity $k$ if and only if $\gcd(q_c^k, f^*f)= q_c^k$. 
\end{crit}

\section{The Jacobian determinant of $t^n$} 

Eilenberg and Niven in \cite{EN}, proved the fundamental theorem of algebra for quaternions, using a degree argument. The key ingredient was Lemma 2 whose proof was depended on the positiveness of $|J(t^n)|$ at the roots of the equation $t^n=\i$.  In this section, we will show that $|J(t^n)|$ is non negative at any point $(x,y,z,w) \in \R^4$.

First, we need a lemma. 

\begin{lem} \label{lem1} Let $u,v,p,q \in \R[x,y,z,w]$ be homogeneous polynomials of the same degree $m \geq 1$ so that $yp=zv$ and $yq=wv$. Let $$F:\R^4 \to \R^4, \qquad 
F(x,y,z,w)=(u,v,p,q)$$ Then, $|J(F)|$ is equal to $y^{-3}mv^2(vu_x-uv_x)$. 
\end{lem}
\proof Since $yp=zv$ and $yq=wv$, we see that $$ \begin{array}{lllll} p_x=\ds\frac zy v_x, \;p_y=\ds\frac{-z}{y^2}v +\ds\frac zy v_y, \;p_z=\ds\frac vy +\ds\frac zy v_z,\; p_w=\ds\frac zy v_w \\ [3ex] q_x=\ds\frac wy v_x, \;q_y=\ds\frac{-w}{y^2}v+\ds\frac wy v_y, \;q_z=\ds\frac wy v_z, \;q_w=\ds\frac vy +\ds\frac wy v_w \end{array}$$
Then, the Jacob{\sf i}an of $F$ is \be \label{JM} J(F)=\left[ \begin{array}{lllll} u_x& u_y&u_z&u_w \\ v_x&v_y&v_z&v_w \\[2ex] \ds\frac zy v_x & \ds\frac{-z}{y^2}v +\ds\frac zy v_y&  \ds\frac vy +\ds\frac zy v_z& \ds\frac zy v_w \\ [3ex] \ds\frac wy v_x& \ds\frac{-w}{y^2}v +\ds\frac wy v_y &\ds\frac wy v_z & \ds\frac vy +\ds\frac wy v_w
\end{array} \right] \ee

Now using the fact that $u,v$ are homogeneous of degree $m$, an easy calculation shows that $|J(F)|=y^{-3}mv^2(vu_x-uv_x)$. \eproof

For $u,v,p,q$ {\sf a's} in Lemma \ref{lem1}, define new polynomials $U,V,P,Q$ by the formula $$ \label{eq4} U +\i V + \j P +\k Q=(x +\i y + \j z +\k w)(u +\i v + \j p +\k q)=(u +\i v + \j p +\k q)(x +\i y + \j z +\k w) $$
Notice that $U,V,P,Q$ are homogeneous of degree $m+1$,  $yP=zV,  yQ=wV$ and $U=xu-yv-zp-wq,  V=yu+xv$. Using the fact that  $yP=zV,  yQ=wV$ we observe that $U$ takes the form $$U=xu -v  \left(\frac{y^2+z^2+w^2}{y} \right )$$

\begin{cor} \label{cor1} Let $U,V,P,Q$ be as above and define $\Phi=(U,V,P,Q)$. Then,  $$ \label{eq5} |J(\Phi)|=\frac{(m+1)V^2[y(x^2+y^2+z^2+w^2)(vu_x-uv_x)+ y^2u^2+v^2(y^2+z^2+w^2)]}{y^4} $$
\end{cor}
\proof We have $$U_x=u+xu_x-v_x \left(\frac{y^2+z^2+w^2}{y} \right ), \quad V_x=yu_x+v+xv_x$$ Now, a calculation shows that $$VU_x-UV_x=\frac{y(x^2+y^2+z^2+w^2)(vu_x-uv_x)+ y^2u^2+v^2(y^2+z^2+w^2)}{y}$$
Thus, the result follows from Lemma \ref {lem1}. \eproof

\begin{cor} \label{cor2} Let $f(t)=t^n=(x +\i y + \j z +\k w)^n=f_1 +\i f_2 +\j f_3 +\k f_4=(f_1, f_2, f_3, f_4)$. Then, $|J(f)| \geq 0$.  
\end{cor} 
\proof For induction purposes, we will rename $f=t^n=f^{[n]}=(f^{[n]_r}), \; r=1,2,3,4$. 
Observe first that $f^{[n]_r}$ are homogeneous of degree $n$ and satisfy the conditions $yf^{[n]_3}=zf^{[n]_2}, yf^{[n]_4}=wf^{[n]_2}$. If $n=1$ we get $|J(t)|=1$, while if $n=2$ we see that $|J(t^2)|=16x^2(x^2+y^2+z^2+w^2) \geq 0$. Let now $n \geq 2$. Then, Corollary \ref{cor1} shows that $|J(t^{n+1})|$ is equal to $$ \frac{(n+1)(f^{[n]_2})^2[y |f^{[1]}|^2  (f^{[n]_2}f^{[n]_1}_x-f^{[n]_1}f^{[n]_2}_x)+ y^2u^2+v^2(y^2+z^2+w^2)]}{y^4}$$
From Lemma \ref{lem1} we see that $y(x^2+y^2+z^2+w^2)(f^{[n]_2}f^{[n]_1}_x-f^{[n]_1}f^{[n]_2}_x) \geq 0$ and thus by induction the result follows. \eproof

\section{The Jacobian of $f(t)$} 

In this section we will compute the Jacobian of a quaternion polynomial $f(t)$. First, we compute the Jacobian of a complex polynomial $a(t) + \i b(t), \, a, b \in \R[t]$ at the complex number $t_0=r+ \i s$ (in the sequel, for a complex polynomial $\phi$, we will denote its Jacobian over $\C$ by $J_{\hbox{\tiny $\C$}}(\phi)$ to differentiate it from the Jacobian of $\phi$ over $\H$). 
Using this, we will  calculate $J(f)(t_0)$. Finally, utilizing the above we compute the Jacobian of $f$ at a point $c_0 \in \H$. 

We begin with the following: 

\begin{lem} \label{lem2} For $n \in \N$, define $a_n=2n-1, \; b_n=2n$ and $c_n=(-1)^{n+1}$. Then,  
\sep {1mm}
\be \label{JM1} J(t^{2n-1})(\i)=\left [ \begin{array}{cccccccc} c_na_n&0&0&0 \\ 0&c_na_n&0&0\\0&0&c_n&0\\0&0&0&c_n \end{array} \right ] \;\; \hbox{and} \;\;  
J(t^{2n})(\i)=\left [ \begin{array}{ccccccc} 0&-c_nb_n&0&0 \\ c_nb_n&0&0&0\\0&0&0&0\\0&0&0&0 \end{array} \right ] \ee
\end{lem}
\proof We will use the notation of Corollary \ref{cor2}. From (\ref{JM}) we see that 
$f_{x_j}^{[n]_3}(\i)=0$ for $x_j=x,y,w$ and $f_{x_j}^{[n]_4}(\i)=0$ for $x_j=x,y,z$. Recall that \be \label{eq8} f^{[k+1]_1}=x f^{[k]} - f^{[k]_2} \left( \frac{y^2+z^2+w^2}{y} \right) \quad \hbox{and} \quad f^{[k+1]_2}=y f^{[k]_1} +x f^{[k]_2} \ee
By using induction and taking into account (\ref{eq8}) we get 
$$ \label{eq9} f^{[2n-1]_1}(\i)=0, \; f^{[2n]_1}(\i)=(-1)^n, \; f^{[2n-1]_2}(\i)=(-1)^{n+1}, \; f^{[2n]_2}(\i)=0 $$
Now $f_x^{[2n+1]_1}=f^{[2n]_1} +x f_x^{[2n]_1} -f_x^{[2n]_2} \left( \frac{y^2+z^2+w^2}{y} \right)$. By induction $f_x^{[2n]_2}(\i)=(-1)^{n+1} (2n)$ and thus $f_x^{[2n+1]_1}(\i)=(-1)^n -(-1)^{n+1} (2n)=(-1)^{n+1+1} (2n+1)$ as required. Similarly, we get $f_y^{[2n-1]_2}(\i)=(-1)^{n+1} (2n-1), \, f_y^{[2n]_1}(\i)=-c_nb_n$ and $f_x^{[2n]_2}(\i)=c_nb_n$. Finally, it is easy to see that $f_{x_j}^{[k]_r}(\i)=0$ when $r=3,4$ and $x_j=z,w$. This finishes the proof. \eproof
Note that if $s \in \R, s \ne 0$, $J(t^k)(\i s)=s^{k-1} J(t^k)(\i)$. 

\begin{cor} \label{cor-complex} Let $g(t)=a(t) + \i b(t), \; a(t), b(t) \in \R[t]$ be a complex polynomial and $r + \i s \in \C-\R$. Then, \be \label{eq-co} J(g)(r + \i s)= \left [ \begin{array}{rrrrrr} \a& -\b &0&0 \\ \b & \a &0&0 \\0&0& \gamma & -\delta \\0&0& \delta & \gamma \end{array} \right ], \quad \hbox{for $\a, \b, \gamma, \delta \in \R$} \ee
\end{cor}
\proof Since $J(g)(r +\i s)=J(g(t+r)(\i s)$, we may assume that $r=0$. We write \be \label{eq-g} g(t)=a_nt^{2n-1} + a_{n-1}t^{2n-3} + \cdots + a_1t + b_mt^{2m} + b_{m-1}t^{2m-2} + \cdots + b_1t^2+ b_0 \ee $a_i, b_j \in \C$. Then, 
$$ \begin{array}{llllll} J(g)(\i s) &=& \sum_{k=1}^n a_k\, J(t^{2k-1})(\i s) + \sum_{l=1}^mb_l \, J(t^{2l})(\i s)\\ [2ex] &=& \sum_{k=1}^n a_k\,s^{2k-2} J(t^{2k-1})(\i) + \sum_{l=1}^m=b_l \, s^{2m-1} J(t^{2l})(\i) \end{array} $$

Now use Lemma \ref{lem2} to get $$ \begin{array}{lllll} \a= \ds\sum_{k=1}^n (-1)^{k+1} s^{2k-2} (2k-1) a_k^1 - \ds\sum_{l=1}^m (-1)^l s^{2l-1} (2l)b_l^2 , \\ [3ex] 
\b= \ds\sum_{k=1}^n (-1)^{k+1} s^{2k-2} (2k-1) a_k^2 + \ds\sum_{l=1}^m (-1)^l s^{2l-1} (2l)b_l^1 , \\ [3ex] \gamma=\ds\sum_{k=1}^n (-1)^{k+1} s^{2k-2} a_k^1\;\;  \hbox{and} \;\; \delta=\ds\sum_{k=1}^n (-1)^{k+1} s^{2k-2} a_k^2
\end{array} $$
where $a_k=a_k^1+ \i a_k^2, \; b_l=b_l^1+ \i b_l^2$. \eproof 

{\sf N}ote that $|J(g)(t) \geq 0$ for all $t \in \C$. On the other hand, $J_{\hbox{\tiny $\C$}}(g)(r + \i s)=\left [ \begin{array}{rrrrrr} \a & -\b \\ \b & \a \end{array} \right]$. Thus, $|J(g)(r +\i s)|=|J_{\hbox{\tiny $\C$}}(g)(r + \i s)| \cdot (\gamma^2 + \delta^2)$. In view of this it may happen that $|J_{\hbox{\tiny $\C$}}(g)(r + \i s)| >0$ while $|J(g)(r + \i s)|=0$. For example, if $g(t)=2t^5+t^3-t+t^2+\i(3t^5+t^3-2t+5)$ and $r=0, s=1$, $\a=-6, \b=-12$ but  $\gamma=\delta=0$. Things, however, become more interesting when $r +\i s$ is a root of $g$. Indeed, we have: 

\begin{rem} \label{rem-complex} For $g, r +\i s$ as in Corollary \ref{cor-complex}, if $|J_{\hbox{\tiny $\C$}}(g)(r +\i s)|>0$ and $|J(g)(r +\i s)|=0$, $g$ is non primitive. 
\end{rem}
\proof Without loss of generality, suppose that $r +\i s=\i$. Write $g$ as in (\ref{eq-g}). Then, $a(\i)=\i \gamma + \sum_{l=1}^m (-1)^lb_l^1 + b_0^1$ and 
$b(\i)=\i \delta + \sum_{l=1}^m (-1)^lb_l^2 + b_0^2$. Since $0=g(\i)=a(\i) +\i b(\i)$ we get $a(\i)=b(\i)=0$. Thus, $a(t), b(t)$ are both divisible by $t^2+1$ which makes $g$ non primitive. \eproof

Now let $f(t)$ be as in (\ref{left}). We write $f(t)=a(t) + \i b(t) +\j ( c(t) -\i d(t))=g(t) +\j h(t) \; a,b,c,d \in \R[t]$. 

\begin{prop} \label{prop-co} For $r +\i s \in \C - \R, f,g,h$ as above, \be \label{eq-quaco} J(f)(r + \i s)= \left [ \begin{array}{rrrrr} \a_1& -\b_1 &-\a_2&\b_2 \\ \b_1 & \a_1 &\b_2&\a_2 \\ \a_3&-\b_3& \a_4 & -\b_4 \\-\b_3&-\a_3& \b_4 & \a_4 \end{array} \right ], \quad \hbox{for $\a_k, \b_k \in \R$} \ee
\end{prop}
\proof We have $J(f)(t)=J(g)(t) + \j J(h)(t)$. Note that $\j J(h)(r + \i s)$ is equal to $$\left [ \begin{array}{rrrrr} 0&0&-1&0\\0&0&0&1\\1&0&0&0 \\0&-1&0&0 \end{array} \right ] \, \left [ \begin{array}{rrrr} \a_3& -\b_3 &0&0 \\ \b_3 & \a_3 &0&0 \\0&0& \a_2 & -\b_2 \\0&0& \b_2 & \a_2 \end{array} \right ]= \left [ \begin{array}{rrrrr} 0&0&-\a_2 & \b_2 \\ 0&0& \b_2 & \a_2 \\ \a_3 & -\b_3 &0&0 \\-\b_3& -\a_3 &0&0 \end{array} \right ]
$$ 
Now the result follows from Corollary \ref{cor-complex}. \eproof 

We finally proceed with the computation of $J(f)$ at $t_0=x_0 +\i y_0 + \j z_0 +\k w_0 \in \H - \C$.  Assuming, as usual, that $x_0=0$, we will ``push'' $t_0$ to $\i s$, with $\i s \sim t_0$ and use  Proposition \ref{prop-co} to find the said Jacobian. 

To achieve the above, for a $c \in \H, |c|=1$, consider the transformation $h_c(t)=ctc*$. This is a rotation in $\R^4$, realized by an orthogonal matrix $A_c$ in the sense that  $h_c(x+ \i y +\j z + \k w)=c(x+ \i y +\j z + \k w)c^*=(h_1 + \i h_2 +\j h_3 + \k h_4)$, then $h(x,y,z,w)=(h_1, h_2, h_3, h_4)=A_c \, (x,y,z,w)^T$. On the other hand, observe that $ctc^*$ commutes with itself and preserves the root structure of $f(t)$. Indeed, if $f(t)$ is as in (\ref{left}), consider a new polynomial $\phi(u)$, 
\be \phi(u)=cf(u)c^*=ca_nc*u^n + ca_{n-1}c^*u^{n-1} + \cdots + ca_1c^*u+ca_0c* \ee 
Note now that if $c_0$ is a root of $f$, $cc_0c^*$ is a root of $\phi$ of the same structure. Divide $f$ by $(t-t_0)$ to get $f(t)=g(t)(t-t_0) + f(t_0)$ and let $\g \in \H$ be so that $\g t_0 \g^*=\i s$. Consider the map $F(x,y,z,w)=(f \circ h_{\g})(x,y,z,w)$. Then, $JF=J(f)(h_{\g}) \cdot J(h_{\g})= 
J(f)(h_{\g}) \cdot A_{\g}$. Therefore, \be \label{eq-qua}  J(f)(t_0)= J(f)(\i s) A_{\g}\ee This along with formula (\ref{eq-quaco}) finishes the task of computing $ J(f)(t_0)$. 

\subsection{Cauchy-Riemann equations}

With the help from the previous section, we are now able to write Cauchy-Riemann (CR) equations for quaternion polynomials. However, the results concerning the CR equations over $\R$ can also be proven independently of the previous ones. The motivation for this is an elementary proof of CR equations for complex polynomials:   

Let $f(t) \in \C[t]$. We write $f(x +\i y)=a(x,y) + \i b(x,y)$.  Then, $a_x=b_y$ and $a_y=-b_x$. 

\proof Let $t_1=(x_1, y_1) \in \C$. Then, $f(t)=f^1(t)(t-t_1) + f(t_1)$. In that case, if $f^1(x +\i y)=a^1(x,y) + \i b^1(x,y)$  we get 
$a(x,y) +\i b(x,y)=(a^1(x,y) +\i b^1(x,y))(x-x_1 +\i (y-y_1))+f(t_1)
=(x-x_1)a^1-(y-y_1)b^1 +\i ((y-y_1)a^1 +(x-x_1)b^1) + f(t_1)$.
Then, $a_x(t_1)=a^1(t_1)$ and $b_y(t_1)=a^1(t_1)$. Similarly, $a_y(t_1)=-b_x(t_1)$.  \eproof

The s{\sf ame} technique can be applied to any monic  quaternion polynomial $f(t)$ for any $t_0 \in \R$.  Indeed, divide $f$ by $t-t_0$ to get $f(t)=g(t)(t-t_0) + f(t_0)=(t^{n-1} + a_{n-2} t^{n-2} + \cdots a_1t+a_0)(t-t_0) + f(t_0)$. Since $t_0$ is real, 
it commutes with any quaternion, and thus we also have $f(t)=(t-t_0)g(t) + f(t_0)$. Thus, if $t=x+\i y +\j z +\k w$ and 
$g(x+\i y +\j z +\k w)=g_1+\i g_2 +\j g_3 + \k g_4$, we see that 
$f(x+\i y +\j z +\k w)=f_1+\i f_2 +\j f_3 + \k f_4=(g_1+\i g_2 +\j g_3 + \k g_4)(x-t_0 +\i y +\j z +\k w) +f(t_0)=(x-t_0 +\i y +\j z +\k w)(g_1+\i g_2 +\j g_3 + \k g_4) +f(t_0)$.  
Now an easy calculation shows that \be \label{jreal} J(f)(t_0)= \left[ \begin{array}{cccccc} g_1(t_0) & -g_2(t_0) &0&0 \\ [1ex] g_2(t_0) & g_1(t_0) &0&0\\0&0&g_1(t_0) &0 \\ [1ex] 0&0&0& g_1(t_0) \end{array} \right ] \ee
with all the other partials equal to zero. This provides the CR equations at a real point $t_0$: 
\be \label{CR-real} \frac{\p f_1}{\p x}=\frac{\p f_2}{\p y}, \;\;\; \frac{\p f_1}{\p y}=-\frac{\p f_2}{\p x}, \;\;\; \frac{\p f_3}{\p z}=\frac{\p f_4}{\p w}=\frac{\p f_1}{\p x} \ee 

The above method, however, fails when $t_0 \in \H - \R$. The reason is that $t$--being a quaternion variable now--does not commute with $t_0$ and thus  $f_1+\i f_2 +\j f_3 + \k f_4 \ne (g_1+\i g_2 +\j g_3 + \k g_4)(x-t_0 +\i y +\j z +\k w) +f(t_0)$. For example, if $f(t)=t^2+ (\i + \j)t-\k=(t+\j)(t+\i)$, there is no linear polynomial 
$g_1+\i g_2 +\j g_3 + \k g_4$ so that $f(x+\i y +\j z +\k w)=(g_1+\i g_2 +\j g_3 + \k g_4)(x +\i(y+1)+\j z +\k w)$.

Next, formula (\ref{eq-co}) provides the CR equations for an $f(t) \in \C[t]$ and $t_0=r +\i s$: \be \label{CR-co} \frac{\p f_1}{\p x}=\frac{\p f_2}{\p y}, \;\;\; \frac{\p f_1}{\p y}=-\frac{\p f_2}{\p x}, \;\;\; \frac{\p f_3}{\p z}=\frac{\p f_4}{\p w}, \;\;\; \frac{\p f_3}{\p w}=-\frac{\p f_4}{\p z}\ee 
Observe that $|J(f)(t_0)| \geq 0$.  

Finally, for the general case, pick $\g \in \H$  so that $\g t_0 \g^*=r+ \i s$.
With the aid of  (\ref{eq-quaco}) and (\ref{eq-qua}) we get the  CR equations for $f$ at $t_0$: 
\be \label{jquat} J(f)(t_0)= \left [ \begin{array}{rrrrr} \a_1& -\b_1 &-\a_2&\b_2 \\ \b_1 & \a_1 &\b_2&\a_2 \\ \a_3&-\b_3& \a_4 & -\b_4 \\-\b_3&-\a_3& \b_4 & \a_4 \end{array} \right ]\cdot A_{\g}, \quad \hbox{where $\a_k, \b_k \in \R$}
\ee

\subsection{$|J(f)|$ at a root of $f$}

In this section we will show that $|J(f)|$ is non negative over $\H$. In particular, we will prove that  if $t_0$ is a root of $f$, $t_0$ is simple if and only if $|J(f)(t_0)| >0$.  Thus, at a multiple root $|J(f)|$ vanishes. Furthermore, we will briefly talk about the local degree of a {\em primitive} $f$ at a root $t_0$ and show how this relates to its multiplicity $\mu(f)(t_0)$. 

Let $t_0=x_0 +\i y_0 + \j z_0 +\k w_0 \in \H$. Divide $f(t)$ by $t-t_0$ to get 
$f(t)=g(t)(t-t_0) + f(t_0)$. After making a suitable transformation, we will assume that $t_0=\i$. Let $g(t)= b_mt^m + \cdots + b_1t +b_0, \; b_k \in \H$. We write $f(t)=b_0(t-\i) + b_1t(t-\i) + \cdots + b_mt^m(t-\i) + f(t_0)$. Let $A$ be the matrix  $$ \label{ma} A= \left [ \begin{array}{rrrrr} 0&-1&0&0\\ 1&0&0&0 \\0&0&0&1\\0&0&-1&0 \end{array} \right ] $$
Notice that $A^2=-I$.  Furthermore, from (\ref{JM1}) we get  $J(t^k(t-\i))(\i)=A^k$ for $k \geq 0$. Therefore, $$ \label{jd} J(f)(\i)=b_0I+b_1A+b_2A^2 + \cdots + b_mA^m=\sum_{k=0}^m (-1)^kb_{2k}I + \sum_
{l=0}^m (-1)^l b_{2l+1}A. 
$$ 
Set $\sum_{k=0}^m (-1)^kb_{2k}=B_e, \; \sum_
{l=0}^m (-1)^l b_{2l+1}=B_o$. We claim that $|B_eI + B_oA| \geq 0$. Indeed, if either of $B_e, B_o$ is zero there is nothing to prove. Suppose then $B_eB_o \ne 0$. Then it is enough to show $|I + CA| \geq 0$ for $C=B_o/B_e$. If $C=a +\i b+ \j c+ \k d,\;$ $I+CA$ takes the form $$\left [ \begin{array}{cccccc} -b+1&-a&d&-c\\a&-b+1&-c&-d\\d&-c&b+1&a\\-c&-d&-a&b+1 \end{array} \right ]$$ 
whose determinant $|I + CA|= q(b^2)=b^4+(-2+2a^2+2c^2+2d^2)b^2+(1-2c^2-2d^2+2a^2+2c^2a^2+2d^2a^2+a^4+d^4+c^4+2d^2c^2)$. The discriminant of $q(b^2)$ is equal to $-16a^2$ and that proves the claim. Moreover, $q(b^2)$ vanishes precisely when $a=0, \; b^2+c^2+d^2=1$. In addition,  notice that if $\d$ is an imaginary unit quaternion, $1+ (a +\i b+ \j c+ \k d)\delta =0$ if and only if $a=0$ and $b^2+c^2+d^2=1$; that is, when $q(b^2)=0$. In short, $|B_eI + B_oA|=0$ if and only $B_e +B_o\,\d=0$, for any imaginary unit quaternion $\d$. 

Let $\g$ be any  imaginary unit quaternion. Then, $g(\g)=B_e + B_o \g$, since $\g^2=-1$. Thus, if $g(\g) \ne 0$, which in turn says that $t=\i$ has multiplicity $1$, $|B_eI + B_oA| > 0$. On the other hand, if $g(\g) =0 $, which means that $\mu(f)(\i) \geq 2$, 
then $J(f)(\g)|$ vanishes, as required. 

Finally, we will consider the relationship between, the local degree of (a primitive) $f$ at its root $c_i$  and its multiplicity $\mu(f)(c_i)$. First, let us recall some facts about degrees of maps. If we compactify $\R^4$ by adding the point $\infty$ to it, then $\R^4 \cup \{\infty\} \equiv S^4$. In that regard, since 
$f(t) \to \infty$ as $t \to \infty$, we get a continuous map $f: S^4 \to S^4$. It is then known that $\deg f=n$, \cite{EN}. Let now $c_1, \cdots, c_m$ be the distinct roots of $f$ and define the local degree $ldf(c_i)$ of $f$ at $c_i$ as follows: 
Choose a ball $B_i$ with center at $c_i$ so that it does not contain any other root of $f$. Let $S_i$ be the boun{\sf da}r{\sf y} of $B_i$. We now can 
define the Gauss map $$G: S_i \to S^3, \quad G(t)=\frac{f(t)}{\|f(t)\|}$$ Then,  $ldf(c_i)=\deg G$. Since the degree of $f$ is $n$, we note that $n=\sum_{i=1}^m ldf(c_i)$. 

If $c_i$ is a simple root, then $ldf(c_i)=1$ since $|J(f)(c_i)| >0$. On the other hand, if $\mu(f)(c_j)=k \geq 2$, we see that $\mu(f^*f)(\a +\i \b)=k$, where $\a+ \i \b \sim c_j$. Since the local degree (over $\C$) of $f*f$ at $\a +\i \b$ is $k$, for a suitable point $\c +\i \d, \; \d \ne 0$ near zero, the equation $f^*f(t)=\c +\i \d$ has $k$ distinct solutions $\g_j +\i \delta_j, \; j=1,2, \cdots, k$ in $\C$, not conjugate to one another. Indeed, if  $\g_1 +\i \delta_1=\g_2 -\i \delta_2$, then 
$f^*f(\g_1+ \i \d_1)=f^*f(\g_2- \i \d_2)=\g+\i \d$. But since $f*f$ is real we should have $\g+\i \d=\g -\i \d$, a contradiction to $\d \ne 0$. Now there exist  {\em distinct} $\zeta_j \sim 
\g_j + \i \d_j$ which are solutions of $f(t)=\a_1 +\i b_1$. Indeed, for if $\zeta_i=\zeta_j$ we have $\eta_1(\g_i +\i \d_i)\eta_1^*=\eta_2(\g_j +\i \d_j)\eta_2^*$ or $\eta_2^*\eta_1(\g_i +\i \d_i)\eta_1^*\eta_2=\g_j +\i \d_j$. Thus   $\eta_2^*\eta_1(\g_i +\i \d_i)\eta_1^*\eta_2 =\g_i + \i \d_i$ which shows $\g_i +\i \d_i=\g_j +\i \d_j$. Therefore, $ldf(c_j)$ must also be $k$. 

\section{Closure} 
In this paper Jacobians of  (left) quaternion polynomials were computed along with their determinants. As a result, Cauchy-Riemman equations were obtained for this type of functions. In addition, local degrees were considered and it was shown that the above commensurates well with the corresponding theory of a (complex) polynomial of a complex variable. Apparently, similar  results can be gotten  for {\em right} quaternion polynomials. 

It is hoped that the fact, of the Jacobian determinant of $f$ at any point in $\H$ being non negative, would help in finding the zeros of $f$ as well local degrees at its zeros. Finally, it would be interesting to investigate whether results of the same nature hold true for general quaternion polynomials of the form $f(t)=a_0ta_1t \cdots ta_n + \phi(t), a_i \in \H, a_i \ne 0$ and $\phi(t)$ being a sum of finite number of similar monomials $b_0tb_1t \cdots tb_k, k <n$.

\def\AMM{{\it Amer.\ Math.\ Monthly\ }}
\def\ACM{{\it Adv.\ Comp.\ Math.\ }}
\def\ACMTMS{{\it ACM Trans.\ Math.\ Software\ }}
\def\ACMTOG{{\it ACM Trans.\ Graphics\ }}
\def\BAMS{{\it Bull.\ Amer.\ Math.\ Soc.\ }}
\def\CAD{{\it Comput.\ Aided Design }}
\def\CAEJ{{\it Comput.\ Aided Eng. J.\ }}
\def\CAGD{{\it Comput.\ Aided Geom.\ Design }}
\def\CAVW{{\it Comput.\ Anim.\ Virt.\ Worlds }}
\def\CG{{\it Computers \& Graphics }}
\def\CVGIP{{\it Comput.\ Vision, Graphics, Image\ Proc.\ }}
\def\GM{{\it Graph.\ Models\ }}
\def\IBMJRD{{\it IBM J.\ Res.\ Develop.\ }}
\def\JCAM{{\it J.\ Comput.\ Appl.\ Math.\ }}
\def\JGG{{\it J.\ Geom.\ Graphics}}
\def\JMAA{{\it J.\ Math.\ Anal.\ Appl.\ }}
\def\JSC{{\it J.\ Symb.\ Comput.\ }}
\def\MC{{\it Math.\ Comp.\ }}
\def\MMAS{{\it Math.\ Methods\ Appl.\ Sci.\ }}
\def\NA{{\it Numer.\ Algor.\ }}
\def\PAMS{{\it Proc.\ Amer.\ Math.\ Soc.\ }}
\def\SIAMJNA{{\it SIAM J.\ Numer.\ Anal.\ }}
\def\SIAMR{{\it SIAM Rev.\ }}
\def\JMAS{{\it J.\ Math.\ Sciences.\ }}
\def\TAMS{{\it Trans.\ Amer.\ Math.\ Soc.\ }}
\def\MJM{{\it Milan\ J.\ Math.\ }}

\begin{flushleft}

\end{flushleft}

\end{document}